\documentclass[a4paper,11pt,reqno]{article}
\usepackage{amsmath, amsfonts, amssymb, amsthm}

\input{epsf.sty}
\usepackage{mathrsfs}
\usepackage{amsmath}
\usepackage{amssymb}
\usepackage{graphicx}
\usepackage{subfigure}
\numberwithin{equation}{section}

\begin{document}

\title{Blow-up phenomena for an integrable two-component Camassa-Holm system with cubic nonlinearity
and peakon solutions}

\author
{Kai $\mbox{Yan}^{1}$ \footnote{E-mail: yankai419@163.com }
\quad
Zhijun $\mbox{Qiao}^{2}$ \footnote{E-mail: qiao@utpa.edu}
\quad and \quad
Yufeng $\mbox{Zhang}^{3}$ \footnote{E-mail: zhangmath@126.com}\\
$^1 \mbox{School}$ of Mathematics and Statistics,\\
Huazhong University of Science and Technology,\\
Wuhan, Hubei 430074, China \\
$^2\mbox{Department}$ of Mathematics, University of Texas-Pan American,\\
Edinburg, Texas 78541, USA\\
$^3 \mbox{College}$ of Sciences, China University of Mining and Technology,\\
Xuzhou, Jiangsu 221116, China}

\date{}
\maketitle

\begin{abstract}
This paper is devoted to an integrable two-component Camassa-Holm system with cubic nonlinearity,
which includes the cubic Camassa-Holm equation (also called the Fokas-Olver-Rosenau-Qiao equation) as a special case.
The one peaked solitons (peakons) and two peakon solutions are described in an explicit formula. Then,
the local well-posedness for the Cauchy problem of the system is studied.
Moreover, we target at the precise blow-up scenario for strong solutions to the system,
and establish a new blow-up result with respect to the initial data.\\


\end{abstract}

\section{Introduction}
The well-known Camassa-Holm (CH) equation \cite{C-H}
\begin{eqnarray*}
m_t+u m_x+2u_x m=0, \quad m=u-u_{xx},
\end{eqnarray*}
has been studied extensively in the past two decades.
It models the unidirectional propagation of shallow water waves over a flat bottom,
and $u(t, x)$ stands for the fluid velocity at time $t$ in the spatial $x$ direction \cite{C-H,D-G-H,J}.
The CH equation is also a model for the propagation of axially symmetric waves in hyperelastic rods \cite{Dai}.
It has a bi-Hamiltonian structure and is completely integrable \cite{C-H,F-F}.
Its solitary waves vanishing at both infinities are peakons \cite{C-H-H},
they are orbitally stable \cite{C-S1}, and the CH equation also possesses the algebro-geometric solutions on a symplectic submanifold \cite{Q-CMP}.
It should be stressed that the peakons replicate a feature that is characteristic for the waves of great height
---waves of largest amplitude that are exact solutions of the governing equations for irrotational water waves,
cf. \cite{Cinvent,C-Eann}.

It has been proved that the Cauchy problem for the CH equation is locally well-posed in \cite{C-Ec,D2}.
Moreover, it has both global strong solutions \cite{Cf,C-Ep,C-Ec}
and solutions with finite time blow-up \cite{Cf,C-Ep,C-Ec,C-E}.
Furthermore, wave breaking (namely, the wave remains bounded while its slope becomes unbounded in finite time \cite{Wh})
is the only way for the CH equation to develop singularities in finite time \cite{C-E}.
On the other hand, it also has global weak solutions after the strong solution to the CH equation
blows up in finite time \cite{B-C2,C-M,X-Z}.
In comparison with the celebrated KdV equation, the advantage of the CH equation lies in the fact that
the CH equation not only has peakon solutions, but also models wave breaking.

The nonlinear term in the CH equation is quadratic.
However, there do exist integrable peakon systems with cubic nonlinearity,
which include the cubic CH equation (also called the Fokas-Olver-Rosenau-Qiao (FORQ) equation)
\begin{eqnarray} \label{cubic-CH}
m_t+((u^2-u^2_x)m)_x=0, \quad m=u-u_{xx},
\end{eqnarray}
and the Novikov equation
\begin{eqnarray*}
m_t+u^2 m_x+3u u_x m=0, \quad m=u-u_{xx}.
\end{eqnarray*}
Eq.(\ref{cubic-CH}) was proposed independently in \cite{Fokas,Fu,OR,Q1}.
Its Lax pair and some explicit soliton solutions have been studied in \cite{Q1}.
Recently, the orbital stability of the train of peakons, as well as the H\"{o}lder continuity
for Eq.(\ref{cubic-CH}) have been studied in \cite{LLQu-2014CH} and \cite{Himonas2014}, respectively.
The Novikov equation has been proposed in \cite{NV1} and its Lax pair, bi-Hamiltonian structure, peakon stability,
well-posedness, blow-up phenomena and global solutions were already studied in \cite{HH,HW,NV1,wu}.

In this paper, we consider the following integrable two-component Camassa-Holm system
with cubic nonlinearity proposed in \cite{SQQ}:
\begin{equation}
\left\{\begin{array}{ll}
m_t+[((uv-u_x v_x)-(uv_x-u_x v))m]_x=0,\\
n_t+[((uv-u_x v_x)-(uv_x-u_x v))n]_x=0,
\end{array}\right.
\end{equation}
where $m=u-u_{xx}$ and $n=v-v_{xx}$. For our convenience, we want to call equation (1.2) the SQQ system.
System (1.2) is a multi-component extension of Eq.(\ref{cubic-CH}),
since it can be reduced to Eq.(\ref{cubic-CH}) as $v \equiv u$.
Integrability of the system (1.2) is shown in \cite{SQQ},
and particularly, this system possesses the following conservation laws:
\begin{eqnarray}\label{conservationlaws}
H_1= \int_{\mathbb{R}} m d x=\int_{\mathbb{R}} u d x,
\quad H_2= \int_{\mathbb{R}} n d x=\int_{\mathbb{R}} v d x,
\end{eqnarray}
$$H_3=\int_{\mathbb{R}} m(v-v_x)dx,$$
and
$$H_4=\int_{\mathbb{R}} (u+u_{x})(v-v_x)^2 m dx.$$
Moreover, some explicit solutions to the system (1.2) such as
the cusped solitons (cuspons) and $W/M$-shape solitons have been given in \cite{SQQ}.
Now, let us set up the Cauchy problem of the SQQ system as follows:
\begin{equation} \label{mn}
\left\{\begin{array}{ll}
m_t+\big((uv-u_x v_x)-(uv_x-u_x v)\big) m_x=-\big((u_x n+ v_x m)+(u n-v m)\big) m,\\
n_t+\big((uv-u_x v_x)-(uv_x-u_x v)\big) n_x=-\big((u_x n+ v_x m)+(u n-v m)\big) n,\\
m(0,x)=m_0(x), \\
n(0,x)=n_0(x).
\end{array}\right.
\end{equation}
By using an approach similar to the one in \cite{Y-Y}, the analytic solutions to the system (\ref{mn})
can readily be proved in both variables, globally in space and locally in time.
However, the main goal of this paper is to study explicit two-peakon solutions to the SQQ system (1.2),
and the blow-up phenomena for strong solutions to the system (\ref{mn}).

To this end, we mainly make good use of the fine structure of the system (\ref{mn}).
When deriving the precise blow-up scenario for the solutions to this system,
we first apply the transport equation theory to obtain a blow-up criterion (see Theorem 4.1 below),
and then exploit the characteristic ODE related to the system (\ref{mn})
to construct some invariant properties of the solutions (see Lemma 4.4 below),
which eventually leads to the precise blow-up scenario (see Theorem 4.2 below).

On the other hand, we directly consider the transport equation in terms of $(u_x n+v_x m)+(u n-v m)$
which is the slope of the term $(uv-u_x v_x)-(uv_x-u_x v)$ (see Lemma 4.5 below)
to derive a new blow-up result with respect to the initial data (see Theorem 4.3 below),
where we observe that the system (\ref{mn}) possesses the conservation laws $H_1$ and $H_2$ in (\ref{conservationlaws}).

The rest of our paper is organized as follows.
In Section 2, we discuss the explicit peakon solutions of the system (1.2).
In Section 3, we state the local well-posedness for the system (\ref{mn}) in Sobolev and Besov spaces.
In Section 4, we derive the precise blow-up scenario and a new blow-up result for the strong solutions to the system (\ref{mn}).

\section{Explicit peakon solutions}
\newtheorem {remark2}{Remark}[section]
\newtheorem{theorem2}{Theorem}[section]
\newtheorem{lemma2}{Lemma}[section]
\newtheorem{corollary2}{Corollary}[section]
\newtheorem{proposition2}{Proposition}[section]

In this section, we derive the explicit peakon solutions and discuss the peakon interactions of the system (1.2).
First of all, taking the convolution with the Green function $p(x)\triangleq \frac{1}{2}e^{-|x|}$
for the Helmholtz operator $(1-\partial_x^2)$, one can rewrite the system (1.2) as the following weak form:
\begin{equation}\label{weak-form}
\left\{\begin{array}{ll}
u_t+(u+u_x)(v-v_x)u_x\\
=-\partial_x p\ast[\left((u+u_x)(v-v_x)\right)_x u_x]- p\ast[\left((u+u_x)(v-v_x)\right)_x u],\\
v_t+(u+u_x)(v-v_x)v_x\\
=-\partial_x p\ast[\left((u+u_x)(v-v_x)\right)_x v_x]- p\ast[\left((u+u_x)(v-v_x)\right)_x v].
\end{array}\right.
\end{equation}

Let us assume the single peakon solution of the system (1.2) shows up in the following form \cite{QXL}:
\begin{eqnarray*}
\left\{
\begin{split}
u(t,x)&=c_1e^{-|x-At|},\\
v(t,x)&=c_2e^{-|x-At|}.
\end{split}
\right.
\label{opteqhuv1}
\end{eqnarray*}
Here $c_1$ and $c_2$ are two arbitrary constants,
and $A$ is the traveling wave speed to be determined by $c_1$ and $c_2$.
The derivatives of the above expressions of $u(t,x)$ and $v(t,x)$ do not exist at $x=At$,
thus they can not directly satisfy the system (1.2) in the classical sense.
However, in the weak sense, we are able to present out the expressions of
$u_t$, $v_t$, $u_x$ and $v_x$ with the help of distribution:
{\small
\begin{eqnarray}
u_t=Ac_1 sgn(x-At)E=-Ac_1E_x, \quad v_t=A c_2 sgn(x-At)E=-Ac_2E_x, \label{uxtm} 
\end{eqnarray}
\begin{eqnarray}
u_x=-c_1sgn(x-At)E=c_1E_x, \quad v_x=-c_2sgn(x-At)E=c_2E_x, \label{vxtn}
\end{eqnarray}
}
where $E=e^{-\mid x-At\mid}$.
Substituting (\ref{uxtm}) and (\ref{vxtn}) into the weak form (\ref{weak-form}), 
taking the integrals with the test function $\phi=\phi(x)$ on $\mathbb{R}$ in the distribution sense, and
noticing the following key identities:
\begin{eqnarray*}
\begin{split}
\int u_t \phi dx&=Ac_1 \int E \phi'(x) dx, \\
\int v_t \phi dx&=Ac_2 \int E \phi'(x) dx,\\
\int(u+u_x)(v-v_x)u_x\phi dx&=0,\\
\int(u+u_x)(v-v_x)v_x\phi dx&=0,\\
\int -\partial_x p\ast[\left((u+u_x)(v-v_x)\right)_x u_x]\phi dx &=2c_1^2c_2 \int p\ast\left(EE_x^2-E_{xx}E_x^2\right)\phi'(x) dx\\
        &=\frac{2c_1^2c_2}{3} \int E\phi'(x) dx,\\
\int -\partial_x p\ast[\left((u+u_x)(v-v_x)\right)_x v_x]\phi dx &=2c_2^2c_1 \int p\ast\left(EE_x^2-E_{xx}E_x^2\right)\phi'(x) dx\\
        &=\frac{2c_2^2c_1}{3} \int E\phi'(x) dx,\\
\int -p\ast[\left((u+u_x)(v-v_x)\right)_x u]\phi dx &=2c_1^2c_2 \int p\ast\left(E^2E_x-EE_xE_{xx}\right)\phi dx=0,\\
\int -p\ast[\left((u+u_x)(v-v_x)\right)_x v]\phi dx &=2c_2^2c_1 \int p\ast\left(E^2E_x-EE_xE_{xx}\right)\phi dx=0,
\end{split}
\end{eqnarray*}
one may immediately
know that $A$ must satisfy
\begin{eqnarray}
A=\frac{2}{3}c_1c_2.
\end{eqnarray}
So, we obtain the following peakon solutions to the system (1.2):
\begin{eqnarray*}
\left\{
\begin{split}
u(t,x)&=c_1e^{-|x-\frac{2}{3}c_1c_2t|},\\
v(t,x)&=c_2e^{-|x-\frac{2}{3}c_1c_2t|}.
\end{split}
\right.
\label{opteqhuv1}
\end{eqnarray*}

Next, we derive the two-peakon solutions of the system (1.2), which possess the following form: 
\begin{eqnarray}
\left\{
\begin{split}
u(t,x)&=p_1e^{-|x-q_1(t)|}+p_2e^{-|x-q_2(t)|},\\
v(t,x)&=r_1e^{-|x-q_1(t)|}+r_2e^{-|x-q_2(t)|},
\end{split}
\right.
\label{opteqhuv1}
\end{eqnarray}
where $p_1, p_2, r_1, r_2, q_1, q_2$ are the functions of $t$ to be determined.

By adopting the procedure similar to the single peakon derivation described above and substituting (\ref{opteqhuv1}) into the weak from (\ref{weak-form}), we have 
\begin{eqnarray}
\left\{\begin{array}{l}
p_{1,t}=p_{2,t}=r_{1,t}=r_{2,t}=0,\\
q_{1,t}=-\frac{1}{3}p_1r_1+\frac{1}{2}\left[p_1r_2\left(sgn(q_1-q_2)-1\right)-p_2r_1\left(sgn(q_1-q_2)+1\right)\right]e^{ -\mid q_1-q_2\mid},\\
q_{2,t}=-\frac{1}{3}p_2r_2+\frac{1}{2}\left[p_1r_2\left(sgn(q_1-q_2)-1\right)-p_2r_1\left(sgn(q_1-q_2)+1\right)\right]e^{ -\mid q_1-q_2\mid}.\\
\end{array}\right. \label{tpSQQ}
\end{eqnarray}
From the first equation of (\ref{tpSQQ}), we know
\begin{eqnarray*}
p_1=A_1, ~~~~ p_2=A_2, ~~~~r_1=B_1,  ~~~~r_2=B_2,
\label{rtp}
\end{eqnarray*}
where $A_1$, $A_2$, $B_1$, and $B_2$ are four integration constants.

If $A_1B_1= A_2B_2$, then we have
\begin{eqnarray*}
\left\{\begin{array}{l}
q_{1}(t)=-2\left\{-\frac{1}{3}A_1B_1+\frac{1}{2}\left[A_1B_2\left(sgn(C_1)-1\right)-A_2B_1\left(sgn(C_1)+1\right)\right]
e^{ -\mid C_1\mid}\right\}t+\frac{C_1}{2},\\
q_{2}(t)=q_{1}(t)-C_1.\\
\end{array}\right. \label{tpSQQc1}
\end{eqnarray*}

\begin{figure}
\begin{minipage}[t]{0.5\linewidth}
\centering
\includegraphics[width=2.2in]{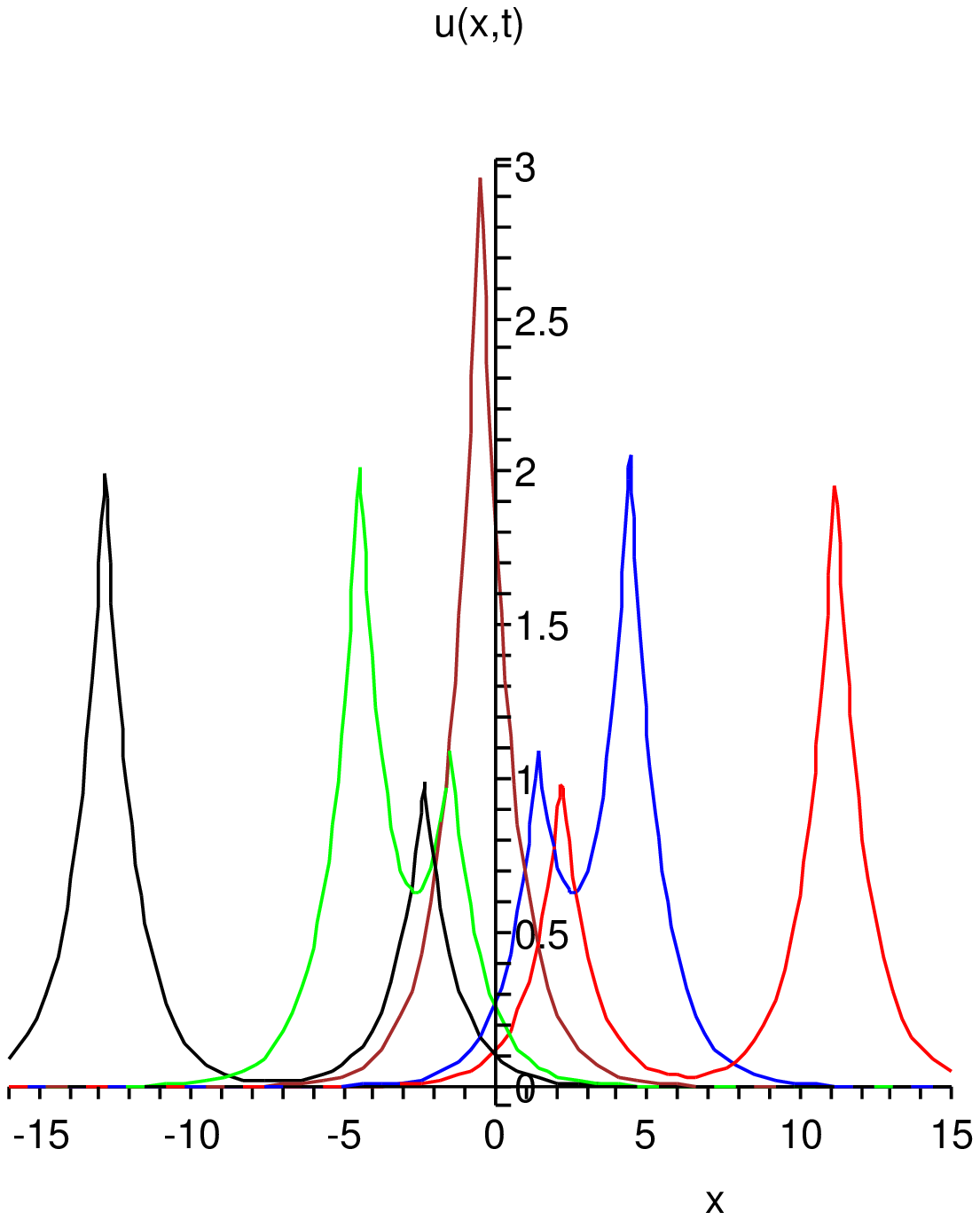}
\caption{\small{ The two-peakon solution for the potential $u(t,x)$ given by (\ref{tpSQQuv}).
Red line: $t=1.5$;
Blue line: $t=0.5$;
Brown line: $t=0$ (collision);
Green line: $t=-0.5$;
Black line: $t=-1.75$. }}
\label{F34}
\end{minipage}
\hspace{2.0ex}
\begin{minipage}[t]{0.5\linewidth}
\centering
\includegraphics[width=2.2in]{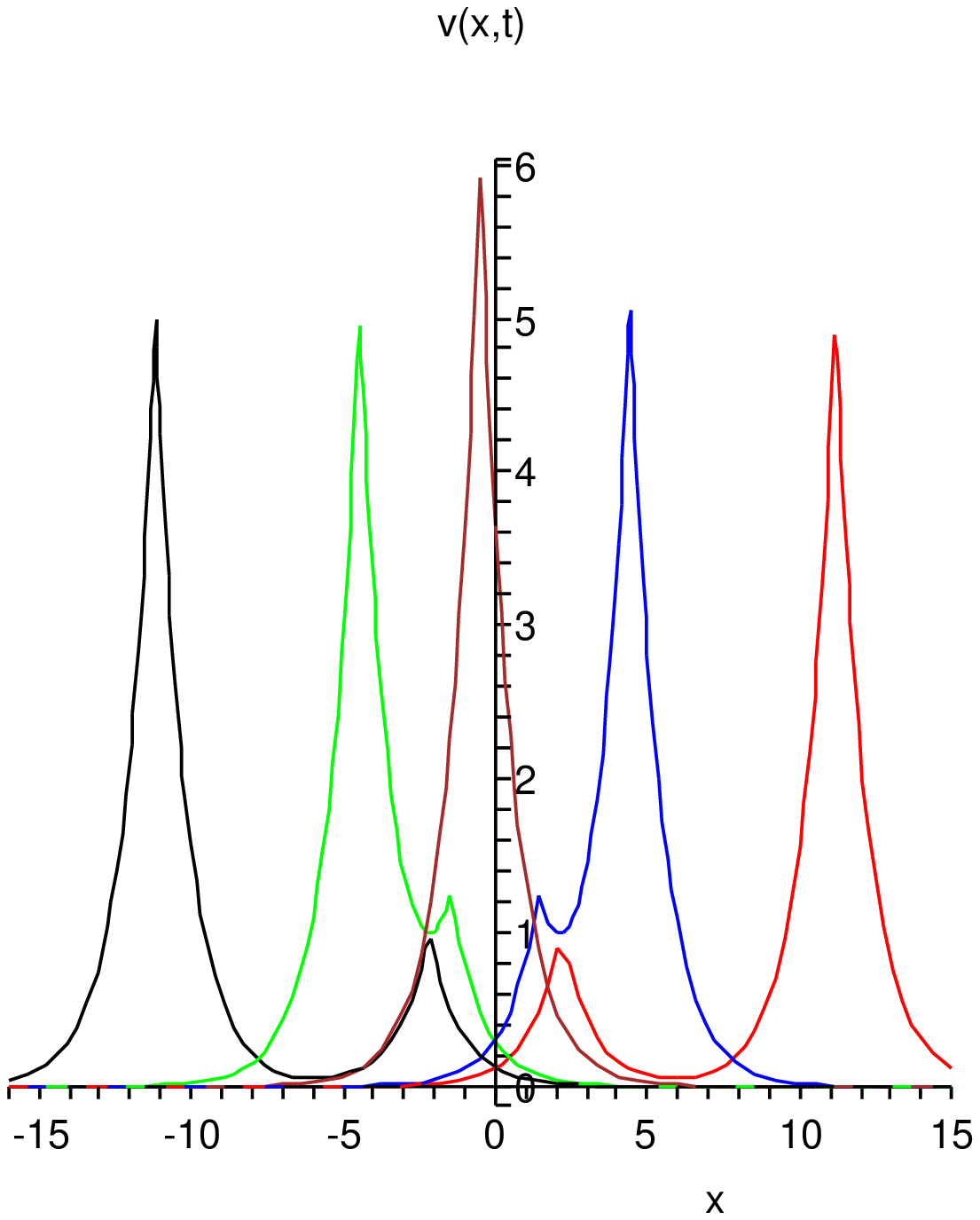}
\caption{\small{ The two-peakon solution for the potential $v(t,x)$ given by (\ref{tpSQQuv}).
Red line: $t=1.5$;
Blue line: $t=0.5$;
Brown line: $t=0$ (collision);
Green line: $t=-0.5$;
Black line: $t=-1.5$. }}
\label{F35}
\end{minipage}
\end{figure}

If $A_1B_1\neq A_2B_2$, then we arrive at
\begin{eqnarray*}
\left\{\begin{array}{l}
q_{1}(t)=\frac{2}{3}A_1B_1t+\Gamma(t),\\
q_{2}(t)=\frac{2}{3}A_2B_2t++\Gamma(t),
\end{array}\right.
\label{tpSQQc2}
\end{eqnarray*}
where
\begin{eqnarray*}
\Gamma(t)
&=&-\frac{3(A_1B_2+A_2B_1)}{2|A_1B_1-A_2B_2|} sgn(t)
\left(e^{-\frac{2}{3}\mid(A_1B_1-A_2B_2)t\mid}-1\right)\\
&&+\frac{3(A_1B_2-A_2B_1)}{2(A_1B_1-A_2B_2)}e^{-\frac{2}{3}\mid(A_1B_1-A_2B_2)t\mid}.
\label{Gamma}
\end{eqnarray*}
In particular, taking $A_1=B_1=1$, $A_2=2$, and $B_2=5$ sends the two-peakon solution to the following form
\begin{eqnarray}
\left\{
\begin{split}
u(t,x)&=e^{-\mid x-q_1(t)\mid}+2e^{-\mid x-q_2(t)\mid},\\
v(t,x)&=e^{-\mid x-q_1(t)\mid}+5e^{-\mid x-q_2(t)\mid},
\end{split}
\right.
\label{tpSQQuv}
\end{eqnarray}
where
\begin{eqnarray*}
\left\{\begin{split}
q_{1}(t)&=\frac{2t}{3}-\frac{7}{6}sgn(t)\left(e^{-6|t|}-1\right)-\frac{1}{2}e^{-6|t|},\\
q_{2}(t)&=\frac{20t}{3}-\frac{7}{6}sgn(t)\left(e^{-6|t|}-1\right)-\frac{1}{2}e^{-6|t|}.
\end{split}
\right.
\label{tpSQQq}
\end{eqnarray*}
For the potential $u(t,x)$, the two-peakon collides at the moment $t=0$, since $q_1(0)=q_2(0)=0$.
For $t>0$, the tall and fast peakon with the amplitude $2$ and peak position $q_2$ chases after
the short and slow peakon with the amplitude $1$ and peak position $q_1$.
At the moment of $t=0$, the two-peakon overlaps.
After the collision ($t<0$), the two-peakon separates, and the tall and fast peakon surpasses the short and slow one.
Similarly, we may discuss the collision of the two-peakon for the potential $v(t,x)$.
See Figures \ref{F34} and \ref{F35} for the two-peakon dynamics of the potentials $u(t,x)$ and $v(t,x)$.

\section{Local well-posedness}
\newtheorem {remark3}{Remark}[section]
\newtheorem{theorem3}{Theorem}[section]
\newtheorem{lemma3}{Lemma}[section]
\newtheorem{corollary3}{Corollary}[section]
\newtheorem{proposition3}{Proposition}[section]

In this section, we will state the local well-posedness for the system (\ref{mn}).
To begin with, applying Kato's semigroup theory \cite{Kato}
and going along the similar line of the proof in \cite{E-Le-Y}, one can readily prove
the following local well-posedness result in Sobolev space.
\begin{theorem3}
Suppose that $(m_0, n_0)\in H^{s}(\mathbb{R})\times H^{s}(\mathbb{R})$ with $s\geq 1$.
There exists a maximal existence time $T=T(||m_0||_{H^{s}(\mathbb{R})}, ||n_0||_{H^{s}(\mathbb{R})})>0$,
and a unique solution $(m, n)$ to the system (\ref{mn}) such that
$$
(m, n)\in C([0,T);H^{s}(\mathbb{R})\times H^{s}(\mathbb{R}))
\cap C^{1}([0,T);H^{s-1}(\mathbb{R})\times H^{s-1}(\mathbb{R})).
$$
Moreover, the solution depends continuously on the initial data,
that is, the mapping $(m_0, n_0)\mapsto (m, n):$
$$ H^{s}(\mathbb{R})\!\times\! H^{s}(\mathbb{R})\!\rightarrow\!
C([0,T); H^{s}(\mathbb{R})\!\times\! H^{s}(\mathbb{R}))\cap
C^{1}([0,T);H^{s-1}(\mathbb{R})\!\times\! H^{s-1}(\mathbb{R}))$$
is continuous.
\end{theorem3}

On the other hand, taking advantage of the transport equation theory, Littlewood-Paley's decomposition
and some fine estimates of Besov spaces \cite{BCD}, one can easily establish the local well-posedness for system (\ref{mn})
in Besov space by means of similar arguments to those in \cite{Yan-Yin-Qiao}.
More precisely,
\begin{theorem3}
Let $1\leq p,r\leq \infty$ and $s>\max (1-\frac{1}{p},\frac{1}{p},\frac 1 2)$ but $s\neq1+\frac{1}{p}$.
Assume that $(m_0,n_0)\in B^s_{p,r}(\mathbb{R})\times B^{s}_{p,r}(\mathbb{R})$,
and set
$$E^s_{p,r}(T)\triangleq C([0,T]; B^s_{p,r}(\mathbb{R}))\cap C^1{([0,T]; B^{s-1}_{p,r}(\mathbb{R}))},
\quad \text{if}\,\ r<\infty,$$
$$E^s_{p,\infty}(T)\triangleq L^{\infty}(0,T; B^s_{p,\infty}(\mathbb{R}))\cap Lip\,(0,T; B^{s-1}_{p,\infty}(\mathbb{R})).$$
Then there exists a time $T>0$ such that $(m,n)\in E^s_{p,r}(T)\times E^{s}_{p,r}(T)$
is the unique solution to the system (\ref{mn}), and the mapping $(m_0,n_0)\mapsto (m,n)$ is continuous from
$B^s_{p,r}(\mathbb{R})\times B^{s}_{p,r}(\mathbb{R})$ into
$$ C([0,T]; B^{s'}_{p,r}(\mathbb{R})\times B^{s'}_{p,r}(\mathbb{R}))\cap
C^1{([0,T]; B^{s'-1}_{p,r}(\mathbb{R})\times B^{s'-1}_{p,r}(\mathbb{R}))}$$
for all $s'<s$ if $r=\infty$, and $s'=s$ if $1\leq r<\infty$.\\
\end{theorem3}

\begin{remark3} \label{remark3-1}
(1) Note that for every $s\in\mathbb{R}$, $B^{s}_{2,2}(\mathbb{R})=H^s(\mathbb{R})$.
Theorem 3.2 holds true in the corresponding Sobolev spaces with $\frac{1} {2}<s\neq\frac{3} {2}$,
which almost improves the result of Theorem 3.1 where $s\geq 1$ is required.
Therefore, the conclusion of Theorem 3.1 holds true for initial data
$(m_0,n_0)\in H^s(\mathbb{R})\times H^s(\mathbb{R})$ with $s>\frac{1} {2}$,
or for all initial data $(u_0,v_0)\in H^s(\mathbb{R})\times H^s(\mathbb{R})$ with $s>\frac{5} {2}$.\\
(2) The above maximal existence time $T$ can be proved independent of the regularity index $s$,
which will be shown in Remark 4.1 below.
\end{remark3}

\section{Blow-up}
\newtheorem {remark4}{Remark}[section]
\newtheorem{theorem4}{Theorem}[section]
\newtheorem{lemma4}{Lemma}[section]
\newtheorem{corollary4}{Corollary}[section]
\newtheorem{proposition4}{Proposition}[section]

In this section, we will derive the precise blow-up scenario of strong solutions to the system (\ref{mn}),
and then state a new blow-up result with respect to the initial data.
Let us first prove a blow-up criterion for the system (\ref{mn}).
For this, we need some a priori estimates of the following transport equation:
\[(TE)\left\{\begin{array}{l}
\partial_t f+v\,\partial_x f=F,\\
f|_{t=0} =f_0.
\end{array}\right.\]

\begin{lemma4} \label{lemma2.1-1}
\cite{BCD}
Let $s>-\frac{1}{2}$. Assume that $f_0\in H^s(\mathbb{R})$, $F\in L^1(0,T; H^s(\mathbb{R}))$, and
$\partial_x v$ belongs to $L^1(0,T; H^{s-1}(\mathbb{R}))$ if $s> \frac{3}{2}$,
or to $L^1(0,T; H^{\frac{1}{2}}(\mathbb{R})\cap L^\infty(\mathbb{R}))$ if $-\frac{1}{2}<s\leq \frac{3}{2}$.
If $f\in L^\infty(0,T; H^s(\mathbb{R}))\cap C([0,T]; \mathcal{S'}(\mathbb{R}))$ solves $(TE)$,
then $f\in C([0,T]; H^s(\mathbb{R}))$.
Moreover, for all $s\neq \frac{3}{2}$, there exists a constant $C=C(s)>0$ such that for all $t\in [0,T]$,
\begin{equation*}
||f(t)||_{H^s}\leq ||f_0||_{H^s}\,+\, \int_0^t ||F(\tau)||_{H^s}d\tau\,+\, C\int_0^t V(\tau)||f(\tau)||_{H^s} d\tau
\end{equation*}
with
\[
V(t)\triangleq \left\{
\begin{array}{ll}
||\partial_x v(t)||_{H^{\frac{1}{2}}\cap L^\infty},& \textrm{if}\quad -\frac{1}{2}<s< \frac{3}{2},\\
||\partial_x v(t)||_{H^{s-1}},& \textrm{if}\quad s> \frac{3}{2}.
\end{array}
\right.
\]
\end{lemma4}

\begin{lemma4} \label{lemma2.2}
\cite{Gui-JFA}
Let $0<s<1$. Assume that $f_0\in H^s(\mathbb{R})$, $F\in L^1(0,T; H^s(\mathbb{R}))$, and
$v,\, \partial_x v \in L^1(0,T; L^\infty(\mathbb{R}))$.
If $f\in L^\infty(0,T; H^s(\mathbb{R}))\cap C([0,T]; \mathcal{S'}(\mathbb{R}))$ solves $(TE)$,
then $f\in C([0,T]; H^s(\mathbb{R}))$.
Moreover, there exists a constant $C=C(s)>0$ such that for all $t\in [0,T]$,
\begin{equation*}
||f(t)||_{H^s}\leq ||f_0||_{H^s}\,+\, \ C \int_0^t ||F(\tau)||_{H^s}d\tau\,+\, C\int_0^t V(\tau)||f(\tau)||_{H^s} d\tau
\end{equation*}
with $V(t)\triangleq ||v(t)||_{L^\infty}+||\partial_x v(t)||_{L^\infty}$.
\end{lemma4}

In addition, the following one-dimensional Morse-type estimates are also required.
\begin{proposition4}\label{proposition3-1}
\cite{BCD,Gui-JFA}
For all $s>0$, there exists a positive constant $C$ independent of $f$ and $g$, such that\\
$$ ||fg||_{H^s(\mathbb{R})}\leq C(||f||_{H^s(\mathbb{R})}||g||_{L^\infty(\mathbb{R})}+||g||_{H^s(\mathbb{R})}||f||_{L^\infty(\mathbb{R})}),$$
and
$$ ||f\partial_x g||_{H^s(\mathbb{R})}\leq
C(||f||_{H^{s+1}(\mathbb{R})}||g||_{L^\infty(\mathbb{R})}+||f||_{L^\infty(\mathbb{R})} ||\partial_x g||_{H^s(\mathbb{R})}).$$
\end{proposition4}

\begin{theorem4}
Let $(m_0,n_0)\in H^s(\mathbb{R})\times H^{s}(\mathbb{R})$ with $s>\frac 1 2$ and
$T$ be the maximal existence time of the solution $(m,n)$ to the system (\ref{mn}), which is guaranteed by Remark \ref{remark3-1}.
If $T<\infty$, then
$$\int_0^T (||m(\tau,\cdot)||_{L^\infty} +||n(\tau,\cdot)||_{L^\infty})^2 d\tau=\infty.$$
\end{theorem4}

\begin{proof}
 We will prove the theorem by induction with respect to the regularity index $s$ $(s>\frac 1 2)$ as follows.

\noindent {\it{Step 1}}. For $s\in(\frac{1}{2},1)$, by Lemma \ref{lemma2.2} and the system (\ref{mn}), we have
\begin{eqnarray*}
&&||m(t)||_{H^{s}}\\
&\leq& ||m_0||_{H^{s}}
+C \int_0^t ||\big((u_x n+v_x m)+(u n-v m)\big)m(\tau)||_{H^s}d\tau
+C\int_0^t ||m(\tau)||_{H^s}\\
&&\times \big(||(uv-u_x v_x)-(uv_x-u_x v)||_{L^{\infty}}+||(u_x n+v_x m)+(u n-v m)||_{L^{\infty}}\big) d\tau
\end{eqnarray*}
and
\begin{eqnarray*}
&&||n(t)||_{H^{s}}\\
&\leq& ||n_0||_{H^{s}}
+C \int_0^t ||\big((u_x n+v_x m)+(u n-v m)\big)n(\tau)||_{H^s}d\tau
+C\int_0^t ||n(\tau)||_{H^s}\\
&&\times \big(||(uv-u_x v_x)-(uv_x-u_x v)||_{L^{\infty}}+||(u_x n+v_x m)+(u n-v m)||_{L^{\infty}}\big) d\tau.
\end{eqnarray*}
Noting that $u=(1-\partial_x^2)^{-1}m=p\ast m$ with $p(x)\triangleq \frac{1}{2}e^{-|x|}\,(x\in\mathbb{R})$, $u_x=(\partial_x p)\ast m$, $u_{xx}=u-m$ and $||p||_{L^1}=||\partial_x p||_{L^1}=1$, together with the Young inequality,
for all $s\in \mathbb{R}$, we have
\begin{eqnarray}
||u||_{L^\infty},\,||u_x||_{L^\infty},\,||u_{xx}||_{L^\infty}\leq C ||m||_{L^\infty}
\end{eqnarray}
and
\begin{eqnarray}
||u||_{H^s},\,||u_x||_{H^s},\,||u_{xx}||_{H^s}\leq C ||m||_{H^s}.
\end{eqnarray}
Similarly, the identity $v=p\ast n$ ensures
\begin{eqnarray}
||v||_{L^\infty},\,||v_x||_{L^\infty},\,||v_{xx}||_{L^\infty}\leq C ||n||_{L^\infty}
\end{eqnarray}
and
\begin{eqnarray}
||v||_{H^s},\,||v_x||_{H^s},\,||v_{xx}||_{H^s}\leq C ||n||_{H^s}.
\end{eqnarray}
Then Proposition \ref{proposition3-1} gives
\begin{eqnarray}
\nonumber &&||\big((u_x n+v_x m)+(u n-v m)\big)m||_{H^s}\\ \nonumber
&\leq& C ||(u_x n+v_x m)+(u n-v m)||_{H^s}||m||_{L^\infty}\\ \nonumber
&&+C||(u_x n+v_x m)+(u n-v m)||_{L^\infty}||m||_{H^s}\\
&\leq& C (||m||_{L^\infty}||n||_{L^\infty}||m||_{H^s}+||m||^2_{L^\infty}||n||_{H^s})
\end{eqnarray}
and
\begin{eqnarray}
\nonumber &&||(uv-u_x v_x)-(uv_x-u_x v)||_{L^{\infty}}+||(u_x n+v_x m)+(u n-v m)||_{L^{\infty}}\\
&&\leq C ||m||_{L^\infty}||n||_{L^\infty}.
\end{eqnarray}
Hence,
\begin{eqnarray*}
||m(t)||_{H^{s}}\leq ||m_0||_{H^{s}}+C \int_0^t
||m(\tau)||_{L^\infty}||n(\tau)||_{L^\infty}||m(\tau)||_{H^s}+||m(\tau)||^2_{L^\infty}||n(\tau)||_{H^s} d\tau.
\end{eqnarray*}
Likewise,
\begin{eqnarray*}
||n(t)||_{H^{s}}\leq ||n_0||_{H^{s}}+C \int_0^t
||n(\tau)||^2_{L^\infty}||m(\tau)||_{H^s}+||m(\tau)||_{L^\infty}||n(\tau)||_{L^\infty}||n(\tau)||_{H^s} d\tau.
\end{eqnarray*}
Thus, we have
\begin{eqnarray}
\,\, && ||m(t)||_{H^{s}}+||n(t)||_{H^{s}}\\ \nonumber
&\leq& ||m_0||_{H^{s}}+||n_0||_{H^{s}}+C \int_0^t
(||m||_{L^\infty}+||n||_{L^\infty})^2 (||m||_{H^s}+||n||_{H^s}) d\tau.
\end{eqnarray}
Taking advantage of Gronwall's inequality, one gets
\begin{eqnarray}
\nonumber && ||m(t)||_{H^{s}}+||n(t)||_{H^{s}}\\
&\leq& (||m_0||_{H^{s}}+||n_0||_{H^{s}}) e^{C \int_0^t (||m||_{L^\infty}+||n||_{L^\infty})^2 d\tau}.
\end{eqnarray}
Therefore, if $T<\infty$ satisfies $\int_0^T (||m(\tau)||_{L^\infty}+||n(\tau)||_{L^\infty})^2 d\tau<\infty$, then we deduce from (4.8) that
\begin{eqnarray}
\limsup\limits_{t\to T}(||m(t)||_{H^s}+||n(t)||_{H^s})<\infty,
\end{eqnarray}
which contradicts the assumption that $T<\infty$ is the maximal existence time. This completes the proof of the theorem for $s\in(\frac 1 2, 1)$.\\

\noindent {\it{Step 2}}. For $s\in[1,\frac{3}{2})$,
applying Lemma \ref{lemma2.1-1} to the first equation of the system (\ref{mn}), we get
\begin{eqnarray*}
||m(t)||_{H^{s}} &\leq& ||m_0||_{H^{s}}
+\int_0^t ||\big((u_x n+v_x m)+(u n-v m)\big)m(\tau)||_{H^s}d\tau\\
&&+C\int_0^t ||m(\tau)||_{H^s} ||(u_x n+v_x m)+(u n-v m)||_{H^{\frac 1 2}\cap L^{\infty}} d\tau
\end{eqnarray*}
Note that
\begin{eqnarray*}
||(u_x n+v_x m)+(u n-v m)||_{H^{\frac 1 2}\cap L^{\infty}}
&\leq& C ||(u_x n+v_x m)+(u n-v m)||_{H^{{\frac {1} {2}}+\varepsilon_0}}\\
&\leq& C ||m||_{H^{{\frac {1} {2}}+\varepsilon_0}} ||n||_{H^{{\frac {1} {2}}+\varepsilon_0}},
\end{eqnarray*}
where $\varepsilon_0\in (0,\frac 1 2)$.
Using (4.5) and the fact that
$H^{\frac 1 2 +\varepsilon_0}(\mathbb{R})\hookrightarrow H^{\frac 1 2}(\mathbb{R})\cap L^\infty(\mathbb{R})$
leads to
\begin{eqnarray*}
||m(t)||_{H^{s}}\leq ||m_0||_{H^{s}}+C \int_0^t
||m||_{H^{{\frac {1} {2}}+\varepsilon_0}} ||n||_{H^{{\frac {1} {2}}+\varepsilon_0}}||m||_{H^s}
+||m||^2_{H^{{\frac {1} {2}}+\varepsilon_0}} ||n||_{H^s} d\tau.
\end{eqnarray*}
For the second equation of the system (\ref{mn}), we can deal with it in a similar way and obtain that
\begin{eqnarray*}
||n(t)||_{H^{s}}\leq ||n_0||_{H^{s}}+C \int_0^t
||n||^2_{H^{{\frac {1} {2}}+\varepsilon_0}} ||m||_{H^s}
+||m||_{H^{{\frac {1} {2}}+\varepsilon_0}} ||n||_{H^{{\frac {1} {2}}+\varepsilon_0}}||n||_{H^s} d\tau.
\end{eqnarray*}
Hence,
\begin{eqnarray*}
\quad  &&||m(t)||_{H^{s}}+||n(t)||_{H^{s}}\\ \nonumber
&\leq& ||m_0||_{H^{s}}+||n_0||_{H^{s}}+C \int_0^t
(||m||_{H^{{\frac {1} {2}}+\varepsilon_0}} +||n||_{H^{{\frac {1} {2}}+\varepsilon_0}})^2
(||m||_{H^s}+||n||_{H^s}) d\tau.
\end{eqnarray*}
Thanks to Gronwall's inequality again, we have
\begin{eqnarray}
&&||m(t)||_{H^{s}}+||n(t)||_{H^{s}}\\ \nonumber
&\leq& (||m_0||_{H^{s}}+||n_0||_{H^{s}})
e^{C \int_0^t (||m||_{H^{{\frac {1} {2}}+\varepsilon_0}} +||n||_{H^{{\frac {1} {2}}+\varepsilon_0}})^2 d\tau}.
\end{eqnarray}
Therefore, if $T<\infty$ satisfies $\int_0^T (||m(\tau)||_{L^\infty}+||n(\tau)||_{L^\infty})^2 d\tau<\infty$,
then we deduce from the uniqueness of the solution to the system (\ref{mn}) and (4.9) with
$\frac{1}{2}+\varepsilon_0\in (\frac{1}{2},1)$ instead of $s$ that
$$||m(t)||_{H^{{\frac {1} {2}}+\varepsilon_0}}+ ||n(t)||_{H^{{\frac {1} {2}}+\varepsilon_0}}\quad
\text{is uniformly bounded in}\quad  t\in (0,T).$$
 This along with (4.10) implies that
\begin{eqnarray}
\limsup\limits_{t\to T}(||m(t)||_{H^s}+||n(t)||_{H^s})<\infty,
\end{eqnarray}
which contradicts the assumption that $T<\infty$ is the maximal existence time. This completes the proof of the theorem for $s\in[1,\frac 3 2)$.\\

\noindent {\it{Step 3}}. For $s\in(1,2)$,
differentiating the system (\ref{mn}) with respect to $x$, we have
\begin{eqnarray*}
&&\partial_t m_x+ \big((uv-u_x v_x)-(uv_x-u_x v)\big)\partial_x m_x \\ \nonumber
&=&-2\big((u_x n+v_x m)+(u n-v m)\big)m_x-\big((u_x n+v_x m)+(u n-v m)\big)_x m\\ \nonumber
&\triangleq& R_1(t,x)
\end{eqnarray*}
and
\begin{eqnarray*}
&&\partial_t n_x+ \big((uv-u_x v_x)-(uv_x-u_x v)\big)\partial_x n_x \\ \nonumber
&=&-2\big((u_x n+v_x m)+(u n-v m)\big)n_x-\big((u_x n+v_x m)+(u n-v m)\big)_x n\\ \nonumber
&\triangleq& R_2(t,x)
\end{eqnarray*}
By Lemma \ref{lemma2.2} with $s-1\in(0,1)$, we get
\begin{eqnarray*}
&&||m_x(t)||_{H^{s-1}}\\
&\leq& ||\partial_x m_0||_{H^{s-1}}
+C \int_0^t ||R_1(\tau)||_{H^{s-1}}d\tau
+C\int_0^t ||m_x(\tau)||_{H^{s-1}}\\
&&\times \big(||(uv-u_x v_x)-(uv_x-u_x v)||_{L^{\infty}}+||(u_x n+v_x m)+(u n-v m)||_{L^{\infty}}\big) d\tau
\end{eqnarray*}
and
\begin{eqnarray*}
&&||n_x(t)||_{H^{s-1}}\\
&\leq& ||\partial_x n_0||_{H^{s-1}}
+C \int_0^t ||R_2(\tau)||_{H^{s-1}}d\tau
+C\int_0^t ||n_x(\tau)||_{H^{s-1}}\\
&&\times \big(||(uv-u_x v_x)-(uv_x-u_x v)||_{L^{\infty}}+||(u_x n+v_x m)+(u n-v m)||_{L^{\infty}}\big) d\tau
\end{eqnarray*}
Due to Proposition \ref{proposition3-1} and (4.1)-(4.4), we have
\begin{eqnarray*}
&&||-2\big((u_x n+v_x m)+(u n-v m)\big)m_x||_{H^{s-1}}\\ \nonumber
&\leq& C\big(||(u_x n+v_x m)+(u n-v m)||_{H^{s}}||m||_{L^{\infty}}\\ \nonumber
&&+||(u_x n+v_x m)+(u n-v m)||_{L^{\infty}}||m_x||_{H^{s-1}}\big)\\ \nonumber
&\leq& C ||m||_{L^{\infty}}||n||_{L^{\infty}}||m||_{H^{s}}
+C ||m||^2_{L^{\infty}}||n||_{H^{s}},
\end{eqnarray*}
and
\begin{eqnarray*}
&&||-\big((u_x n+v_x m)+(u n-v m)\big)_x m||_{H^{s-1}}\\ \nonumber
&\leq& C\big(||m||_{H^{s}}||(u_x n+v_x m)+(u n-v m)||_{L^{\infty}}\\ \nonumber
&&+||m||_{L^{\infty}}||(u_x n+v_x m)+(u n-v m)||_{H^{s}}\big)\\ \nonumber
&\leq& C ||m||_{L^{\infty}}||n||_{L^{\infty}}||m||_{H^{s}}
+C ||m||^2_{L^{\infty}}||n||_{H^{s}},
\end{eqnarray*}
which together with (4.6) yields
\begin{eqnarray*}
||m_x(t)||_{H^{s-1}}\leq ||m_0||_{H^{s}}+C \int_0^t
||m||_{L^\infty}||n||_{L^\infty}||m||_{H^s}+||m||^2_{L^{\infty}}||n||_{H^{s}} d\tau.
\end{eqnarray*}
Likewise,
\begin{eqnarray*}
||n_x(t)||_{H^{s-1}}\leq ||n_0||_{H^{s}}+C \int_0^t
||n||^2_{L^{\infty}}||m||_{H^{s}}+||m||_{L^\infty}||n||_{L^\infty}||n||_{H^s} d\tau.
\end{eqnarray*}
Thus, we have
\begin{eqnarray*}
\,\, && ||m_x(t)||_{H^{s-1}}+||n_x(t)||_{H^{s-1}}\\ \nonumber
&\leq& ||m_0||_{H^{s}}+||n_0||_{H^{s}}+C \int_0^t
(||m||_{L^\infty}+||n||_{L^\infty})^2 (||m||_{H^s}+||n||_{H^s}) d\tau.
\end{eqnarray*}
This along with (4.7) with $s-1\in (0,1)$ instead of $s$ ensures
\begin{eqnarray*}
\,\, && ||m(t)||_{H^{s}}+||n(t)||_{H^{s}}\\ \nonumber
&\leq& ||m_0||_{H^{s}}+||n_0||_{H^{s}}+C \int_0^t
(||m||_{L^\infty}+||n||_{L^\infty})^2(||m||_{H^s}+||n||_{H^s}) d\tau.
\end{eqnarray*}
Similar to Step 1, we can easily prove the theorem for $s\in(1,2)$.\\

\noindent {\it{Step 4}}. For $s=k\in\mathbb{N}$ and $k\geq2$,
differentiating the system (\ref{mn}) $k-1$ times with respect to $x$, we get
\begin{eqnarray*}
&&[\partial_t+\big((uv-u_x v_x)-(uv_x-u_x v)\big)\partial_x]\partial^{k-1}_x m\\
&=&-\sum\limits_{l=0}^{k-2} C^l_{k-1}\partial^{k-l-1}_x\big((uv-u_x v_x)-(uv_x-u_x v)\big)\partial^{l+1}_x m\\
&&-\partial^{k-1}_x [\big((u_x n+v_x m)+(u n- v m)\big)m]\\ \nonumber
&\triangleq& F_1(t,x)
\end{eqnarray*}
and
\begin{eqnarray*}
&&[\partial_t+\big((uv-u_x v_x)-(uv_x-u_x v)\big)\partial_x]\partial^{k-1}_x n\\
&=&-\sum\limits_{l=0}^{k-2} C^l_{k-1}\partial^{k-l-1}_x\big((uv-u_x v_x)-(uv_x-u_x v)\big)\partial^{l+1}_x n\\
&&-\partial^{k-1}_x [\big((u_x n+v_x m)+(u n- v m)\big)n]\\ \nonumber
&\triangleq& F_2(t,x),
\end{eqnarray*}
which together with Lemma \ref{lemma2.1-1} imply
\begin{eqnarray*}
||\partial^{k-1}_x m(t)||_{H^1} &\leq& ||m_0||_{H^{k}}
+\int_0^t ||F_1(\tau)||_{H^1}d\tau+C\int_0^t ||m(\tau)||_{H^k}\\
&&\times ||(u_x n+v_x m)+(u n-v m)||_{H^{\frac 1 2}\cap L^{\infty}} d\tau
\end{eqnarray*}
and
\begin{eqnarray*}
||\partial^{k-1}_x n(t)||_{H^1} &\leq& ||n_0||_{H^{k}}
+\int_0^t ||F_2(\tau)||_{H^1}d\tau+C\int_0^t ||n(\tau)||_{H^k}\\
&&\times ||(u_x n+v_x m)+(u n-v m)||_{H^{\frac 1 2}\cap L^{\infty}} d\tau.
\end{eqnarray*}
Making use of  Proposition \ref{proposition3-1} and (4.1)-(4.4) again, one infers
\begin{eqnarray}
 &&||-\sum\limits_{l=0}^{k-2} C^l_{k-1}\partial^{k-l-1}_x\big((uv-u_x v_x)-(uv_x-u_x v)\big)
\partial^{l+1}_x m||_{H^{1}}
\\ \nonumber
&\leq& C(k)\sum\limits_{l=0}^{k-2}
(||\partial^{k-l-1}_x\big((uv-u_x v_x)-(uv_x-u_x v)\big)||_{L^{\infty}}||m||_{H^{l+2}}\\ \nonumber
&&+||\partial^{k-l-1}_x\big((uv-u_x v_x)-(uv_x-u_x v)\big)||_{H^{1}}||\partial^{l+1}_x m||_{L^{\infty}})\\ \nonumber
&\leq& C(k)\sum\limits_{l=0}^{k-2}
(||(uv-u_x v_x)-(uv_x-u_x v)||_{H^{k-l-\frac{1}{2}+\varepsilon_0}}||m||_{H^{l+2}}\\ \nonumber
&&+||(uv-u_x v_x)-(uv_x-u_x v)||_{H^{k-l}}||m||_{H^{l+\frac{3}{2}+\varepsilon_0}})\\ \nonumber
&\leq& C(k) (||(uv-u_x v_x)-(uv_x-u_x v)||_{H^{k-\frac{1}{2}+\varepsilon_0}}||m||_{H^{k}} \\ \nonumber
&&+||(uv-u_x v_x)-(uv_x-u_x v)||_{H^{k}}||m||_{H^{k}-\frac{1}{2}+\varepsilon_0})\\ \nonumber
&\leq& C(k) (||m||_{H^{k-\frac{1}{2}+\varepsilon_0}}||n||_{H^{k-\frac{1}{2}+\varepsilon_0}}||m||_{H^{k}}
+||m||^2_{H^{k-\frac{1}{2}+\varepsilon_0}}||n||_{H^{k}}),
\end{eqnarray}
\begin{eqnarray}
\nonumber &&||-\partial^{k-1}_x [\big((u_x n+v_x m)+(u n- v m)\big)m]||_{H^{1}}\\ \nonumber
&\leq& C||\big((u_x n+v_x m)+(u n- v m)\big)m||_{H^{k}}\\ \nonumber
&\leq& C||m||_{L^{\infty}}||n||_{L^{\infty}}||m||_{H^{k}}
+||(u_x n+v_x m)+(u n- v m)||_{H^{k}}||m||_{L^{\infty}}\\
&\leq& C(||m||_{H^{k-\frac{1}{2}+\varepsilon_0}}||n||_{H^{k-\frac{1}{2}+\varepsilon_0}}||m||_{H^{k}}
+||m||^2_{H^{k-\frac{1}{2}+\varepsilon_0}}||n||_{H^{k}}),
\end{eqnarray}
and
\begin{eqnarray*}
&&||(u_x n+v_x m)+(u n-v m)||_{H^{\frac 1 2}\cap L^{\infty}}\\
&\leq& C ||(u_x n+v_x m)+(u n-v m)||_{H^{k-\frac{1}{2}+\varepsilon_0}}\\
&\leq& C ||m||_{H^{k-\frac{1}{2}+\varepsilon_0}}||n||_{H^{k-\frac{1}{2}+\varepsilon_0}},
\end{eqnarray*}
where $\varepsilon_0\in (0,\frac 1 2)$ and we used the fact that
\begin{eqnarray}
H^{k-\frac 1 2 +\varepsilon_0}(\mathbb{R})\hookrightarrow H^{\frac 1 2 +\varepsilon_0}(\mathbb{R})
\hookrightarrow H^{\frac 1 2}(\mathbb{R})\cap L^\infty(\mathbb{R})\quad \text{with}\quad  k\geq 2.
\end{eqnarray}
Thus, we get
\begin{eqnarray*}
||\partial^{k-1}_x m(t)||_{H^{1}}
&\leq& ||m_0||_{H^{k}}+C \int_0^t
(||m||_{H^{k-\frac{1}{2}+\varepsilon_0}}||n||_{H^{k-\frac{1}{2}+\varepsilon_0}}||m||_{H^k}\\
&&+||m||^2_{H^{k-\frac{1}{2}+\varepsilon_0}}||n||_{H^{k}}) d\tau.
\end{eqnarray*}
Similarly,
\begin{eqnarray*}
||\partial^{k-1}_x n(t)||_{H^{1}}
&\leq& ||n_0||_{H^{k}}+C \int_0^t
(||n||^2_{H^{k-\frac{1}{2}+\varepsilon_0}}||m||_{H^{k}}\\
&&+||m||_{H^{k-\frac{1}{2}+\varepsilon_0}}||n||_{H^{k-\frac{1}{2}+\varepsilon_0}}||n||_{H^k}) d\tau.
\end{eqnarray*}
Then,
\begin{eqnarray*}
\,\, && ||\partial^{k-1}_x m(t)||_{H^{1}}+||\partial^{k-1}_x n(t)||_{H^{1}}\\ \nonumber
&\leq& ||m_0||_{H^{k}}+||n_0||_{H^{k}}+C \int_0^t
(||m||_{H^{k-\frac{1}{2}+\varepsilon_0}}+||n||_{H^{k-\frac{1}{2}+\varepsilon_0}})^2(||m||_{H^k}+||n||_{H^k}) d\tau,
\end{eqnarray*}
which together with Gronwall's inequality and (4.10) with $s=1$ imply
\begin{eqnarray}
\nonumber &&||m(t)||_{H^{k}}+||n(t)||_{H^{k}}\\
&\leq& (||m_0||_{H^{k}}+||n_0||_{H^{k}})
e^{C \int_0^t (||m||_{H^{k-\frac{1}{2}+\varepsilon_0}}+||n||_{H^{k-\frac{1}{2}+\varepsilon_0}})^2 d\tau}.
\end{eqnarray}
If $T<\infty$ satisfies $\int_0^T (||m(\tau)||_{L^\infty}+||n(\tau)||_{L^\infty})^2 d\tau<\infty$,
applying Step 3 with $\frac 3 2 +\varepsilon_0\in(1,2)$ and by induction with respect to $k\geq 2$,
we see that $||m(t)||_{H^{k-\frac{1}{2}+\varepsilon_0}}+||n(t)||_{H^{k-\frac{1}{2}+\varepsilon_0}}$
is uniformly bounded in $t\in (0,T)$. By (4.15), we have
\begin{eqnarray}
\limsup\limits_{t\to T}(||m(t)||_{H^k}+||n(t)||_{H^k})<\infty,
\end{eqnarray}
which contradicts the assumption that $T<\infty$ is the maximal existence time. This completes the proof of the theorem for $s=k\in\mathbb{N}$ and $k\geq 2$.\\

\noindent {\it{Step 5}}. For $s\in(k,k+1),\,k\in\mathbb{N} $ and $k\geq 2$,
differentiating the system (\ref{mn}) $k$ times with respect to $x$, we get
\begin{eqnarray*}
&&[\partial_t+\big((uv-u_x v_x)-(uv_x-u_x v)\big)\partial_x]\partial^{k}_x m\\
&=&-\sum\limits_{l=0}^{k-1} C^l_{k}\partial^{k-l}_x\big((uv-u_x v_x)-(uv_x-u_x v)\big)\partial^{l+1}_x m\\
&&-\partial^{k}_x [\big((u_x n+v_x m)+(u n-v m)\big)m]\\ \nonumber
&\triangleq& G_1(t,x)
\end{eqnarray*}
and
\begin{eqnarray*}
&&[\partial_t+\big((uv-u_x v_x)-(uv_x-u_x v)\big)\partial_x]\partial^{k}_x n\\
&=&-\sum\limits_{l=0}^{k-1} C^l_{k}\partial^{k-l}_x\big((uv-u_x v_x)-(uv_x-u_x v)\big)\partial^{l+1}_x n\\
&&-\partial^{k}_x [\big((u_x n+v_x m)+(u n-v m)\big)n]\\ \nonumber
&\triangleq& G_2(t,x),
\end{eqnarray*}
which together with Lemma \ref{lemma2.2} with $s-k\in(0,1)$ imply
\begin{eqnarray*}
&&||\partial^{k}_x m(t)||_{H^{s-k}}\\
&\leq& ||\partial^{k}_x m_0||_{H^{s-k}}
+ C \int_0^t ||G_1(\tau)||_{H^{s-k}}d\tau
+C \int_0^t ||\partial^{k}_x m(\tau)||_{H^{s-k}}\\
&&\times \big(||(uv-u_x v_x)-(uv_x-u_x v)||_{L^{\infty}}+||(u_x n+v_x m)+(u n-v m)||_{L^{\infty}}\big) d\tau
\end{eqnarray*}
and
\begin{eqnarray*}
&&||\partial^{k}_x n(t)||_{H^{s-k}}\\
&\leq& ||\partial^{k}_x n_0||_{H^{s-k}}
+ C \int_0^t ||G_2(\tau)||_{H^{s-k}}d\tau
+C \int_0^t ||\partial^{k}_x n(\tau)||_{H^{s-k}}\\
&&\times \big(||(uv-u_x v_x)-(uv_x-u_x v)||_{L^{\infty}}+||(u_x n+v_x m)+(u n-v m)||_{L^{\infty}}\big) d\tau.
\end{eqnarray*}
By (4.14) and using the procedure similar to (4.12)-(4.13), we have
\begin{eqnarray*}
&&||-\sum\limits_{l=1}^{k-1} C^l_{k}\partial^{k-l}_x\big((uv-u_x v_x)-(uv_x-u_x v)\big)
\partial^{l+1}_x m||_{H^{s-k}}\\
&&+||-\partial^{k}_x [\big((u_x n+v_x m)+(u n-v m)\big)m]||_{H^{s-k}}\\ \nonumber
&\leq& C(k)(||m||_{H^{k-\frac{1}{2}+\varepsilon_0}}||n||_{H^{k-\frac{1}{2}+\varepsilon_0}}||m||_{H^{s}}
+||m||^2_{H^{k-\frac{1}{2}+\varepsilon_0}}||n||_{H^{s}}),
\end{eqnarray*}
and
\begin{eqnarray*}
&&||- C^0_{k}\partial^{k}_x\big((uv-u_x v_x)-(uv_x-u_x v)\big) m_x||_{H^{s-k}}\\ \nonumber
&\leq& C (||m_x||_{H^{s-k+1}}||\partial^{k-1}_x\big((uv-u_x v_x)-(uv_x-u_x v)\big)||_{L^\infty}\\
&&+||m_x||_{L^\infty}||\partial^{k}_x\big((uv-u_x v_x)-(uv_x-u_x v)\big)||_{H^{s-k}})\\ \nonumber
&\leq& C (||m||_{H^{s-k+2}}||(uv-u_x v_x)-(uv_x-u_x v)||_{H^{k-\frac{1}{2}+\varepsilon_0}}\\
&&+||m||_{H^{k-\frac{1}{2}+\varepsilon_0}}||(uv-u_x v_x)-(uv_x-u_x v)||_{H^{s}})\\ \nonumber
&\leq& C (||m||_{H^{k-\frac{1}{2}+\varepsilon_0}}||n||_{H^{k-\frac{1}{2}+\varepsilon_0}}||m||_{H^{s}}
+||m||^2_{H^{k-\frac{1}{2}+\varepsilon_0}}||n||_{H^{s}}).
\end{eqnarray*}
Thus, we obtain
\begin{eqnarray*}
\,\, && ||\partial^{k}_x m(t)||_{H^{s-k}}+||\partial^{k}_x n(t)||_{H^{s-k}}\\ \nonumber
&\leq& ||m_0||_{H^{s}}+||n_0||_{H^{s}}+C \int_0^t
(||m||_{H^{k-\frac{1}{2}+\varepsilon_0}}+||n||_{H^{k-\frac{1}{2}+\varepsilon_0}})^2(||m||_{H^s}+||n||_{H^s}) d\tau.
\end{eqnarray*}
This along with (4.7) with $s-k\in (0,1)$ instead of $s$ lead to
\begin{eqnarray*}
\,\, && ||m(t)||_{H^{s}}+|| n(t)||_{H^{s}}\\ \nonumber
&\leq& ||m_0||_{H^{s}}+||n_0||_{H^{s}}+C \int_0^t
(||m||_{H^{k-\frac{1}{2}+\varepsilon_0}}+||n||_{H^{k-\frac{1}{2}+\varepsilon_0}})^2(||m||_{H^s}+||n||_{H^s}) d\tau.
\end{eqnarray*}
By using Gronwall's inequality, Step 3 with $\frac 3 2+\varepsilon_0\in (1,2)$
and the similar argument as shown in Step 4, we can arrive at the desired result.

In summary, the above 5 steps complete the proof of the theorem.
\end{proof}

\begin{remark4}
The maximal existence time $T$ in Theorem 4.1 can be chosen independent of the regularity index $s$.
Indeed, let $(m_0,n_0)\in H^{s}\times H^{s}$ with $s>\frac{1}{2}$ and some $s'\in (\frac 1 2,s)$.
Then Remark 3.1 ensures that there exists a unique $H^s\times H^s$ (resp., $H^{s'}\times H^{s'}$) solution
$(m_s,n_s)$ (resp., $(m_{s'},n_{s'})$) to the system (\ref{mn}) with the maximal existence time $T_s$ (resp., $T_{s'}$).
Since $H^s\hookrightarrow H^{s'}$, it follows from the uniqueness that
$T_s\leq T_{s'}$ and $(m_s,n_s)\equiv (m_{s'},n_{s'})$ on $[0, T_s)$.
On the other hand, if we suppose that $T_s < T_{s'}$, then $(m_{s'},n_{s'})\in C([0,T_s]; H^{s'}\times H^{s'})$.
Hence $(m_{s},n_{s})\in L^2(0,T_s; L^{\infty}\times L^{\infty})$, which is a contradiction to Theorem 4.1.
Therefore, $T_s = T_{s'}$.
\end{remark4}

Now we turn our attention to the precise blow-up scenario for sufficiently regular solutions to the system (\ref{mn}).
For this, motivated by \cite{Cf,McKean}, we first consider the characteristic ordinary differential equation as follows:
\begin{equation}
\left\{\begin{array}{ll}
\frac{d q(t,x)}{dt}=\big((u v-u_x v_x)-(u v_x-u_x v)\big)(t,q(t,x)),\quad\quad &(t,x)\in(0,T)\times \mathbb{R}, \\
q(0,x)=x,\quad\quad  &x\in\mathbb{R},
\end{array}\right.
\end{equation}
for the flow generated by $(u v-u_x v_x)-(u v_x-u_x v)$.

The following lemmas are very crucial to study the blow-up phenomena of strong solutions to the system (\ref{mn}).
\begin{lemma4}
Let $(m_0,n_0)\in H^s(\mathbb{R})\times H^s(\mathbb{R})$ with $s> \frac 1 2$
and $T>0$ be the maximal existence time of the corresponding solution $(m,n)$ to the system (\ref{mn}).
Then Eq.(4.17) has a unique solution $q\in C^1([0,T)\times \mathbb{R};\mathbb{R})$.
Moreover, the mapping $q(t,\cdot)$ is an increasing diffeomorphism of $\mathbb{R}$ with
\begin{eqnarray}
q_{x}(t,x)=\exp\left(\int_0^t \big((u_x n+v_x m)+(u n-v m)\big)(s,q(s,x))ds\right)>0,
\end{eqnarray}
for all $(t,x)\in [0,T)\times \mathbb{R}$.
\end{lemma4}

\begin{proof}
Since $(u,v)\in C([0,T); H^s(\mathbb{R})\times H^s(\mathbb{R}))
\cap C^1([0,T); H^{s-1}(\mathbb{R})\times H^{s-1}(\mathbb{R}))$ with $s>\frac 5 2$, it follows from the fact
$H^{s-1}(\mathbb{R})\hookrightarrow Lip(\mathbb{R})$ with $s>\frac 5 2$ that $(u v-u_x v_x)-(u v_x-u_x v)$ is bounded and Lipschitz continuous in the space variable $x$ and of class $C^1$ in time variable $t$. Then the classical ODE theory ensures that Eq.(4.17) has a unique solution $q\in C^1([0,T)\times \mathbb{R};\mathbb{R})$.\\
Differentiating Eq.(4.17) with respect to $x$ gives
\begin{equation*}
\left\{\begin{array}{ll}
\frac{d q_x(t,x)}{dt}=\big((u_x n+v_x m)+(u n-v m)\big)(t,q(t,x))q_x(t,x),\quad\quad
&(t,x)\in(0,T)\times \mathbb{R}, \\
q_x(0,x)=1,\quad\quad  &x\in\mathbb{R},
\end{array}\right.
\end{equation*}
which leads to (4.18).\\
On the other hand, $\forall\, t<T$, by the Sobolev embedding theorem, we have
$$\sup\limits_{(s,x)\in [0,T)\times \mathbb{R}} \left|\big((u_x n+v_x m)+(u n-v m)\big)(s,x)\right| < \infty,$$
which along with (4.18) implies that there exists a constant $C>0$ such that
$$q_x(t,x)\geq e^{-Ct},\quad \forall\, (t,x)\in [0,T)\times \mathbb{R}.$$
This implies that  the mapping $q(t,\cdot)$ is an increasing diffeomorphism of $\mathbb{R}$ before blow-up.
Therefore, we complete the proof of Lemma 4.3.
\end{proof}

\begin{lemma4}
Let $(m_0,n_0)\in H^s(\mathbb{R})\times H^s(\mathbb{R})$ with $s> \frac 1 2$
and $T>0$ be the maximal existence time of the corresponding solution $(m,n)$ to the system (\ref{mn}).
Then we have
\begin{equation}
m(t,q(t,x))q_x(t,x)=m_0(x),
\end{equation}
and
\begin{equation}
n(t,q(t,x))q_x(t,x)=n_0(x).
\end{equation}
for all $(t,x)\in [0,T)\times \mathbb{R}$.\\
Moreover, if there exists a $C>0$ such that for all $(t,x)\in [0,T)\times \mathbb{R}$,
$$\big((u_x n+v_x m)+(u n-v m)\big)(t,x)\geq -C,$$
then for all $t\in [0,T)$,
\begin{equation}
||m(t,\cdot)||_{L^\infty}\leq C e^{Ct}||m_0||_{H^s} \quad \text{and} \quad
||n(t,\cdot)||_{L^\infty}\leq C e^{Ct}||n_0||_{H^s}.
\end{equation}
\end{lemma4}

\begin{proof}
Differentiating the left-hand side of (4.19)-(4.20) with
respect to $t$ and making use of (4.17)-(4.18) and the system (\ref{mn}), we have
\begin{eqnarray*}
&&\frac{d}{dt}{(m(t,q(t,x))q_x(t,x))}\\
&=&(m_t(t,q)+m_x(t,q)q_t(t,x))q_{x}(t,x)+m(t,q)q_{xt}(t,x)\\
&=&[m_t+\big((u v-u_x v_x)-(u v_x-u_x v)\big) m_x+\big((u_x n+v_x m)\\
&&+(u n-v m)\big) m](t,q(t,x)) q_{x}(t,x)\\
&=& 0
\end{eqnarray*}
and
\begin{eqnarray*}
&&\frac{d}{dt}{(n(t,q(t,x))q_x(t,x))}\\
&=&(n_t(t,q)+n_x(t,q)q_t(t,x))q_{x}(t,x)+n(t,q)q_{xt}(t,x)\\
&=&[n_t+\big((u v-u_x v_x)-(u v_x-u_x v)\big) n_x+\big((u_x n+v_x m)\\
&&+(u n-v m)\big) n](t,q(t,x)) q_{x}(t,x)\\
&=& 0,
\end{eqnarray*}
which proves (4.19) and (4.20).
By Lemma 4.3, in view of (4.18)-(4.20), the assumption of the lemma, and the fact
$H^s(\mathbb{R})\hookrightarrow L^{\infty}(\mathbb{R})$ as $s>\frac 1 2$,
we obtain for all $t\in[0,T)$,
\begin{eqnarray*}
||m(t,\cdot)||_{L^{\infty}}=||m(t,q(t,\cdot))||_{L^{\infty}}=||q_x^{-1}(t,\cdot) m_0(\cdot)||_{L^{\infty}}
\leq C e^{Ct}||m_0||_{H^{s}}
\end{eqnarray*}
and
\begin{eqnarray*}
||n(t,\cdot)||_{L^{\infty}}=||n(t,q(t,\cdot))||_{L^{\infty}}=||q_x^{-1}(t,\cdot) n_0(\cdot)||_{L^{\infty}}
\leq C e^{Ct}||n_0||_{H^{s}},
\end{eqnarray*}
which complete the proof of the lemma.
\end{proof}

The following theorem shows the precise blow-up scenario for sufficiently regular solutions to the system (\ref{mn}).
\begin{theorem4}
Let $(m_0,n_0)\in H^s(\mathbb{R})\times H^s(\mathbb{R})$ with $s> \frac 1 2$
and $T>0$ be the maximal existence time of the corresponding solution $(m,n)$ to the system (\ref{mn}).
Then the solution $(m,n)$ blows up in finite time
if and only if
$$\liminf_{t\rightarrow T}\inf_{x\in\mathbb{R}}\{\left((u_x n+ v_x m)+(u n-v m)\right)(t,x)\}=-\infty.$$
\end{theorem4}

\begin{proof}
Assume that the solution $(m,n)$ blows up in finite time ($T<\infty$) and
there exists a $C>0$ such that
$$\big((u_x n+v_x m)+(u n-v m)\big)(t,x)\geq -C,\quad \forall\, (t,x)\in [0,T)\times \mathbb{R}.$$
By (4.21), we have
$$\int_0^T (||m(t)||_{L^\infty}+||n(t)||_{L^\infty})^2 d t
\leq C^2 T e^{2CT}(||m_0||_{H^{s}}+||n_0||_{H^{s}})^2<\infty,$$
which contradicts to Theorem 4.1.

On the other hand, by (4.1)-(4.4) and the Sobolev embedding theorem, we can see that if
$$\liminf_{t\rightarrow T}\inf_{x\in\mathbb{R}}\{\left((u_x n+ v_x m)+(u n-v m)\right)(t,x)\}=-\infty,$$
then the solution $(m,n)$ must blow up in finite time. This completes the proof of the theorem.
\end{proof}

\begin{remark4}
If $v\equiv u$, then Theorem 4.2 covers the corresponding result in \cite{Gui-CMP}.
\end{remark4}

In order to have a new blow-up criterion with respect to the initial data of strong solutions to the system (\ref{mn}),
we directly investigate the transport equation in terms of $(u_x n+ v_x m)+(u n-v m)$ which is the slope of
$(u v-u_x v_x)-(u v_x-u_x v)$.

\begin{lemma4}
Let $(m_0,n_0)\in H^s(\mathbb{R})\times H^s(\mathbb{R})$ with $s> \frac 1 2$
and $T>0$ be the maximal existence time of the corresponding solution $(m,n)$ to the system (\ref{mn}).
Set $M=M(t,x)\triangleq \big((u_x n+ v_x m)+(u n-v m)\big)(t,x)$.
Then for all $(t,x)\in [0,T)\times \mathbb{R}$, we have
\begin{eqnarray}
&& M_t+\big((u v-u_x v_x)-(u v_x-u_x v)\big) M_x\\ \nonumber
&=&- M^2-n(1-\partial^2_x)^{-1}\big((u_x+u) M\big)
-m(1-\partial^2_x)^{-1}\big((v_x-v) M\big)\\ \nonumber
&&-n\partial_x(1-\partial^2_x)^{-1}\big((u_x+u) M\big)
+m\partial_x(1-\partial^2_x)^{-1}\big((v_x-v) M\big).
\end{eqnarray}
Moreover, if we assume that $(m_0, n_0)\in L^1(\mathbb{R})\times L^1(\mathbb{R})$ and
$m_0(x),\, n_0(x)\geq 0$ for all $x\in\mathbb{R}$, then
\begin{eqnarray}
|u_x(t,x)| \leq u(t,x) \leq  H_1= \int_{\mathbb{R}} m_0(x) d x,
\end{eqnarray}
\begin{eqnarray}
|v_x(t,x)| \leq v(t,x) \leq  H_2= \int_{\mathbb{R}} n_0(x) d x,
\end{eqnarray}
and
\begin{eqnarray}
M_t+\big((u v-u_x v_x)-(u v_x-u_x v)\big) M_x \leq - M^2+\delta (m+n),
\end{eqnarray}
for all $(t,x)\in [0,T)\times \mathbb{R}$,
where $\delta \triangleq \frac 3 2 (H_1+H_2)^3$ is a positive constant and $H_1, H_2$ are two conservation laws in (\ref{conservationlaws}).
\end{lemma4}

\begin{proof}
As per Remark 4.1, we here may assume $s\geq 3$ to prove the lemma. Firstly, we have
\begin{eqnarray}
&& M_t+\big((u v-u_x v_x)-(u v_x-u_x v)\big) M_x \\ \nonumber
&=&u_{xt}n+v_{xt}m+u_t n-v_t m+u_x n_t+v_x m_t+u n_t-v m_t
+\big((u v-u_x v_x)\\ \nonumber
&&-(u v_x-u_x v)\big) \times (u_x n_x+v_x m_x+u_{xx}n+v_{xx}m+u_x n-v_x m+u n_x-v m_x).
\end{eqnarray}
From the system (\ref{mn}), we infer that
\begin{eqnarray*}
&&(1-\partial^2_x)[u_t+\big((u v-u_x v_x)-(u v_x-u_x v)\big) u_x]\\
&=& m_t+(1-\partial^2_x)[\big((u v-u_x v_x)-(u v_x-u_x v)\big) u_x]\\
&=& -\big((u v-u_x v_x)-(u v_x-u_x v)\big) m_x- M m+\big((u v-u_x v_x)-(u v_x-u_x v)\big)u_x\\
&&-\partial^2_x [\big((u v-u_x v_x)-(u v_x-u_x v)\big) u_x]\\
&=& -\big((u v-u_x v_x)-(u v_x-u_x v)\big)(m_x-u_x)- M m\\
&&-\partial_x [M u_x+\big((u v-u_x v_x)-(u v_x-u_x v)\big) (u-m)]\\
&=& -2 u M+ M m-u_x M_x\\
&=& -u M-(u_x M)_x.
\end{eqnarray*}
Hence,
\begin{eqnarray}
u_t+\big((u v-u_x v_x)-(u v_x-u_x v)\big) u_x = -(1-\partial^2_x)^{-1}\big(u M+(u_x M)_x\big).
\end{eqnarray}
Likewise,
\begin{eqnarray}
v_t+\big((u v-u_x v_x)-(u v_x-u_x v)\big) v_x = -(1-\partial^2_x)^{-1}\big(v M+(v_x M)_x\big).
\end{eqnarray}
By virtue of (4.27)-(4.28) and the system (\ref{mn}), we have
\begin{eqnarray}
&& u_{xt}n+v_{xt}m \\ \nonumber
&=&-\big((u v-u_x v_x)-(u v_x-u_x v)\big)(u_{xx}n+v_{xx}m)
-n\partial_x(1-\partial^2_x)^{-1}(u M)\\ \nonumber
&&-m\partial_x(1-\partial^2_x)^{-1}(v M)
-n(1-\partial^2_x)^{-1}(u_x M)-m(1-\partial^2_x)^{-1}(v_x M),
\end{eqnarray}
\begin{eqnarray}
&& u_{t}n-v_{t}m \\ \nonumber
&=&-\big((u v-u_x v_x)-(u v_x-u_x v)\big)(u_{x}n-v_{x}m)
-n\partial_x(1-\partial^2_x)^{-1}(u_x M)\\ \nonumber
&&+m\partial_x(1-\partial^2_x)^{-1}(v_x M)
-n(1-\partial^2_x)^{-1}(u M)+m(1-\partial^2_x)^{-1}(v M),
\end{eqnarray}
\begin{eqnarray}
&&u_x n_t+v_x m_t\\ \nonumber
&=&-\big((u v-u_x v_x)-(u v_x-u_x v)\big)(u_x n_x+v_x m_x)-(u_x n+v_x m)M,
\end{eqnarray}
and
\begin{eqnarray}
&& u n_t-v m_t\\ \nonumber
&=&-\big((u v-u_x v_x)-(u v_x-u_x v)\big)(u n_x-v m_x)-(u n-v m)M.
\end{eqnarray}
Applying (4.26) and (4.29)-(4.32), one gets (4.22).

Since $m_0(x),\, n_0(x)\geq 0$ for all $x\in\mathbb{R}$, it follows from (4.18)-(4.20) that
\begin{eqnarray}
m(t,x),\, n(t,x)\geq 0, \quad \forall \, (t,x)\in[0,T)\times \mathbb{R}.
\end{eqnarray}
Note that
$$u(t,x)=(1-\partial^2_x)^{-1}m(t,x)=(p\ast m)(t,x)=\frac 1 2\int_{\mathbb{R}} e^{-|x-y|}m(t,y)d y.$$
Then
$$u(t,x)=\frac{e^{-x}}{2}\int^x_{-\infty} e^y m(t,y)d y+\frac{e^{x}}{2}\int_x^{\infty} e^{-y} m(t,y)d y$$
and
$$u_x(t,x)=-\frac{e^{-x}}{2}\int^x_{-\infty} e^y m(t,y)d y+\frac{e^{x}}{2}\int_x^{\infty} e^{-y} m(t,y)d y,$$
which together with (4.33) yields
\begin{eqnarray}
0\leq u(t,x)+u_x(t,x)= \int_x^{\infty} e^{x-y} m(t,y)d y\leq H_1 =\int_{\mathbb{R}} m_0(x) d x
\end{eqnarray}
and
\begin{eqnarray*}
0\leq u(t,x)-u_x(t,x)= \int^x_{-\infty} e^{y-x} m(t,y)d y\leq H_1 =\int_{\mathbb{R}} m_0(x) d x.
\end{eqnarray*}
Hence, we prove (4.23).\\
Similarly, in view of $v(t,x)=(1-\partial^2_x)^{-1} n(t,x)=(p\ast n) (t,x)$, one gets
\begin{eqnarray*}
0\leq v(t,x)+v_x(t,x)= \int_x^{\infty} e^{x-y} n(t,y)d y\leq H_2 =\int_{\mathbb{R}} n_0(x) d x
\end{eqnarray*}
and
\begin{eqnarray}
0\leq v(t,x)-v_x(t,x)= \int^x_{-\infty} e^{y-x} n(t,y)d y\leq H_2 =\int_{\mathbb{R}} n_0(x) d x,
\end{eqnarray}
which implies (4.24).

Noting that
\begin{eqnarray*}
|\partial_x (1-\partial^2_x)^{-1}f(x)|
&=& \left|\frac 1 2 \int_{\mathbb{R}}sgn (x-y) e^{-|x-y|} f(y)d y\right|\\
&\leq& \frac 1 2 \int_{\mathbb{R}} e^{-|x-y|} |f(y)|d y\\
&=& (p\ast |f|)(x),
\end{eqnarray*}
and applying (4.23)-(4.24), (4.34)-(4.35), together with the facts $u=p\ast m,\, v=p\ast n$ again, one infers
\begin{eqnarray*}
&& -n(1-\partial^2_x)^{-1}\big((u_x+u) M\big)-m(1-\partial^2_x)^{-1}\big((v_x-v) M\big)\\
&\leq& n \big(p\ast ((u_x+u)(v-v_x)m)\big)+m \big(p\ast ((u_x+u)(v-v_x)n)\big)\\
&\leq& H_1 H_2 \big(n(p\ast m)+m(p\ast n)\big)\\
&\leq& H_1 H_2 (H_1 n +H_2 m)\\
&\leq& H_1 H_2 (H_1+ H_2)(m+n),
\end{eqnarray*}
and
\begin{eqnarray*}
&& -n\partial_x(1-\partial^2_x)^{-1}\big((u_x+u) M\big)+m\partial_x(1-\partial^2_x)^{-1}\big((v_x-v) M\big)\\
&\leq& n \big(p\ast |(u_x+u)M|\big)+m \big(p\ast |(v-v_x)M|\big)\\
&\leq& H_1^2 n(p\ast n)+H_1 H_2 n(p\ast m)+H_2^2 m(p\ast m)+H_1 H_2 m(p\ast n)\\
&\leq& 2 H_1 H_2 (H_1+ H_2)(m+n),
\end{eqnarray*}
which along with (4.22) imply (4.25). Therefore, we prove the lemma.
\end{proof}

With Lemma 4.5 in hand, a new blow-up result with regard to the initial data follows to conclude this section.

\begin{theorem4}
Suppose that $(m_0,n_0)\in \big(H^s(\mathbb{R})\cap L^1(\mathbb{R})\big)\times
\big(H^s(\mathbb{R})\cap L^1(\mathbb{R})\big)$ with $s> \frac 1 2$
and $T>0$ be the maximal existence time of the corresponding solution $(m,n)$ to the system (\ref{mn}).
Assume that $m_0(x),\, n_0(x)\geq 0$ for all $x\in\mathbb{R}$,
and $m_0(x_0),\, n_0(x_0)> 0$ for some $x_0\in\mathbb{R}$.
Set $M(t)\triangleq M(t,q(t,x_0))$ and $N(t)\triangleq (m+n)(t,q(t,x_0))$.
If
\begin{eqnarray}
M(0)<-\sqrt{2 \delta N(0)},
\end{eqnarray}
where $\delta$ is the positive constant defined in Lemma 4.5,
then the solution $(m, n)$ blows up at a time $T_0\in (0,T_1]$ with
\begin{eqnarray*}
T_1\triangleq -\frac{M(0)+\sqrt{M^2(0)-2\delta N(0)}}{\delta N(0)}.
\end{eqnarray*}
Moreover, when $T_0=T_1$, we have the estimate of the blow-up rate as follows
\begin{eqnarray}
\liminf\limits_{t\rightarrow T_0^-} \big((T_0-t) M(t)\big)\leq -1.
\end{eqnarray}
\end{theorem4}

\begin{proof}
In view of Remark 4.1, we here may assume $s\geq 3$ to prove the theorem.
By (4.17), (4.25) and the assumption of the theorem, we have
\begin{eqnarray}
\nonumber \frac{d}{dt} M(t)&=&\frac{d}{dt} M(t,q(t,x_0))\\ \nonumber
&=&\left(M_t+\big((u v-u_x v_x)-(u v_x-u_x v)\big)M_x\right)(t,q(t,x_0))\\
&\leq& -M^2(t)+\delta N(t).
\end{eqnarray}
From the system (\ref{mn}) and Eq.(4.17), we get
\begin{eqnarray}
\frac{d}{dt} N(t)&=&\frac{d}{dt} m(t,q(t,x_0))+\frac{d}{dt} n(t,q(t,x_0))\\ \nonumber
&=& -M(t)N(t).
\end{eqnarray}
Apparently, (4.18)-(4.20) and the assumption imply $N(t)>0$ for all $t\in[0,T)$.
By (4.38) and (4.39), we have
\begin{eqnarray*}
N(t)\frac{d}{dt} M(t)-M(t)\frac{d}{dt} N(t)\leq \delta N^2(t),
\end{eqnarray*}
which gives
$\frac{d}{dt} \left(\frac{M(t)}{N(t)}\right) \leq \delta.$
Integrating from $0$ to $t$ yields
\begin{eqnarray*}
\frac{M(t)}{N(t)}\leq \frac{M(0)}{N(0)}+ \delta t,
\end{eqnarray*}
or hence
\begin{eqnarray}
M(t)\leq \left(\frac{M(0)}{N(0)}+\delta t\right)N(t),
\end{eqnarray}
which along with (4.39) leads to
\begin{eqnarray}
\frac{d}{dt}\left(\frac{1}{N(t)}\right)\leq \frac{M(0)}{N(0)}+ \delta t.
\end{eqnarray}
Integrating from $0$ to $t$ gives
\begin{eqnarray}
\nonumber 0<\frac{1}{N(t)}
&\leq& \frac{1}{N(0)}+\int_0^t \left(\frac{M(0)}{N(0)}+ \delta \tau \right) d\tau\\ \nonumber
&=&\frac{\delta}{2} \left(t^2+\frac{2 M(0)}{\delta N(0)}t+\frac{2}{\delta N(0)}\right)\\
&=&\frac{\delta}{2}(t-T_1) (t-T_2),
\end{eqnarray}
with
\begin{eqnarray}
0< T_1<-\frac{M(0)}{\delta N(0)}=\frac{T_1+T_2}{2}< T_2,
\end{eqnarray}
and $T_1, T_2$ are two different roots of the equation
$t^2+\frac{2 M(0)}{\delta N(0)}t+\frac{2}{\delta N(0)}=0$, which is ensured by (4.36).
Thanks to (4.42), one can find some $T_0\in (0, T_1]$ such that
$$N(t)\rightarrow +\infty,\quad \text{as}\,t\rightarrow T_0.$$
In view of (4.40) and (4.43), one deduces
$$\inf\limits_{x\in\mathbb{R}} \big((u_x n+v_x m)+(u n-v m)\big)(t,x)
\leq M(t)\rightarrow -\infty,
\quad \text{as}\,t\rightarrow T_0.$$
So, according to Theorem 4.2, the solution $(m,n)$ blows up at the time $T_0\in(0,T_1]$.

On the other hand, when $T_0=T_1$, thanks to (4.40) and (4.42)-(4.43),
one infers that for all $0<t<T_0$,
\begin{eqnarray*}
(T_0-t)M(t) &\leq& (T_0-t)\left(\frac{M(0)}{N(0)}+ \delta t\right) N(t)\\
&\leq& (T_0-t)\left(\frac{M(0)}{N(0)}+ \delta t\right) \frac{2}{\delta (t-T_1) (t-T_2)}\\
&=&-2\frac {t-\frac{T_1+T_2}{2}}{t-T_2},
\end{eqnarray*}
which yields (4.37).
Therefore, we complete the proof of the theorem.
\end{proof}

\bigskip
\noindent\textbf{Acknowledgments}
The authors thank the referees for their valuable comments and suggestions.
This work was supported the National Natural Science Foundation of China (Grant Nos. 11171295 and 61328103).
The author (Z. Qiao) also thanks the U.S. Department of Education GAANN project (P200A120256) to support UTPA mathematics graduate program.

\end{document}